\documentclass[aos,MSNbibl,seceqn,dvips]{arximspdf}

\doi{10.1214/12-AOS1020} 
\referstodoi{10.1214/11-AOS949}
\volume{40}
\issue{4}
\pubyear{2012}
\firstpage{2005}
\lastpage{2013}

\makeatletter
\newcommand{\tr}{\operatorname{tr}}
\newcommand{\lh}{\ell}
\makeatother

\begin{document}
\begin{frontmatter}

\title{Rejoinder: Latent variable graphical model selection via convex
optimization}
\runtitle{Rejoinder}

\begin{aug}
\author[A]{\fnms{Venkat} \snm{Chandrasekaran}\corref{}\ead[label=e1]{venkatc@caltech.edu}},
\author[B]{\fnms{Pablo A.} \snm{Parrilo} \ead[label=e2]{parrilo@mit.edu}}
\and
\author[B]{\fnms{Alan~S.}~\snm{Willsky}\ead[label=e3]{willsky@mit.edu}}
\runauthor{V. Chandrasekaran, P. A. Parrilo and A. S. Willsky}
\affiliation{California Institute of Technology, Massachusetts Institute of Technology and Massachusetts Institute of Technology}
\address[A]{V. Chandrasekaran\\
Department of Computing\\
\quad and Mathematical Sciences\\
California Institute of Technology\\
Pasadena, California 91106\\
USA\\
\printead{e1}}
\address[B]{P. A. Parrilo\\
A.~S. Willsky\\
Laboratory for Information\\
\quad and Decision Systems\\
Department of Electrical Engineering\\
\quad and Computer Science\\
Massachusetts Institute of Technology\\
Cambridge, Massachusetts 02139\\
USA\\
\printead{e2}\\
\phantom{E-mail:}
\printead*{e3}} 
\end{aug}

\received{\smonth{5} \syear{2012}}



\end{frontmatter}
%

\section{Introduction}\label{sec:intro}

We thank all the discussants for their careful reading of our paper,
and for their insightful critiques. We would also like to thank the
editors for organizing this discussion. Our paper contributes to the
area of high-dimensional statistics which has received much attention
over the past several years across the statistics, machine learning and
signal processing communities. In this rejoinder we clarify and
comment on some of the points raised in the discussions. Finally, we
also remark on some interesting challenges that lie ahead in latent
variable modeling.

Briefly, we considered the problem of latent variable graphical model
selection in the Gaussian setting. Specifically, let $X$ be a
zero-mean Gaussian random vector in $\mathbb{R}^{p+h}$ with $O$ and $H$
representing disjoint subsets of indices in $\{1,\ldots,p+h\}$ with $|O|
= p$ and $|H| = h$. Here the subvector $X_O$ represents the observed
variables and the subvector $X_H$ represents the latent variables.
Given samples of only the variables $X_O$, is it possible to
consistently perform model selection? We noted that if the number of
latent variables $h$ is small relative to $p$ and if the conditional
statistics of the observed variables $X_O$ conditioned on the latent
variables $X_H$ are given by a sparse graphical model, then the
marginal concentration matrix of the observed variables $X_O$ is given
as the sum of a sparse matrix and a low-rank matrix. As a first step
we investigated the identifiability of latent variable Gaussian
graphical models---effectively, this question boils down to one of
uniquely decomposing the sum of a sparse matrix and a low-rank matrix
into the individual components. By studying the geometric properties
of the algebraic varieties of sparse and low-rank matrices, we provided
natural sufficient conditions for identifiability and gave statistical
interpretations of these conditions. Second, we proposed the following
regularized maximum-likelihood estimator to decompose the concentration
matrix into sparse and low-rank components:
%
\begin{eqnarray}\label{eq:sdp}
&&(\hat{S}_n,\hat{L}_n) = \arg\min_{S,L}    -\lh(S-L; \Sigma ^n_O) +
  \lambda_n   \bigl(\gamma\|S\|_{1} + \tr(L)\bigr) \nonumber\\ [-8pt]\\ [-8pt]
  &&\hphantom{(\hat{S}_n,\hat{L}_n) =}\quad\mbox{s.t. }    S-L \succ0,    L \succeq0.\nonumber
\end{eqnarray}
Here $\Sigma_O^n$ represents the sample covariance formed from $n$
samples of the observed variables, $\ell$ is the Gaussian
log-likelihood function, $\hat{S}_n$ represents the estimate of the
conditional graphical model of the observed variables conditioned on
the latent variables, and $\hat{L}_n$ represents the extra correlations
induced due to marginalization over the latent variables. The $\ell_1$
norm penalty induces sparsity in $\hat{S}_n$ and the trace norm penalty
induces low-rank structure in $\hat{L}_n$. An important feature of this
estimator is that it is given by a convex program that can be solved
efficiently. Our final contribution was to establish the
high-dimensional consistency of this estimator under suitable
assumptions on the Fisher information underlying the true model (in the
same spirit as irrepresentability conditions for sparse model selection
\cite{RavWRY2008,ZhaY2006}).

\section{Alternative estimators}

A number of the commentaries described alternative formulations for
estimators in the latent variable setting.

\subsection{EM-based methods}


The discussions by Yuan and by Lauritzen and Meinshausen describe an
EM-based alternative in which the rank of the matrix $L$ is explicitly
constrained:
%
\begin{eqnarray}\label{eq:em}
(\hat{S}_n,\hat{L}_n) = \arg\min_{S,L}      -\lh(S-L; \Sigma ^n_O) +
  \lambda_n   \|S\|_{1} \nonumber\\ [-8pt]\\ [-8pt]
  \eqntext{\mbox{s.t. }      S-L \succ0,    L \succeq0,
   \operatorname{rank}(L) \leq\mathrm{r}.}
\end{eqnarray}
The experimental results based on this approach seem quite promising,
and certainly deserve further investigation. On the one hand, we
should reiterate that the principal motivation for our convex
optimization based formulation was to develop a method for latent
variable modeling with provable statistical \emph{and} computational
guarantees. One of the main drawbacks of EM-based methods is the
existence of local optima in the associated variational formulations,
thus leading to potentially different solutions depending on the
initial point. On the other hand, one of the reasons for the positive
empirical behavior observed by Yuan and by Lauritzen and Meinshausen
may be that all the local optima in the experimental settings
considered by the authors may be ``good'' models. Such behavior has in
fact been rigorously characterized recently for certain nonconvex
estimators in some missing data problems \cite{LohW2011}.

One of the motivations for the EM proposal of Yuan and of Lauritzen and
Meinshausen seems to be that there are fairly mature and efficient
solvers for the graphical lasso.\vadjust{\goodbreak} As our estimator is relatively newer
and as its properties are better understood going forward, we expect
that more efficient solvers will be developed for (\ref{eq:sdp}) as
well. Indeed, the LogdetPPA solver \cite{WanST2010} that we cite in
our paper already scales to instances involving several hundred
variables, while more recent efforts \cite{MaXZ2012} have resulted in
algorithms that scale to instances with several thousand variables.

\subsection{Thresholding estimators}


Ren and Zhou propose and analyze an interesting thresholding based
estimator for decomposing a concentration matrix into sparse and
low-rank components. They apply a two-step procedure---$\ell_1$ norm
thresholding followed by trace norm thresholding---to obtain the
sparse component followed by the low-rank component. Roughly speaking,
this two-step estimator can be viewed as the application of the first
cycle of a block coordinate descent procedure to compute our estimator
that alternately updates the sparse and low-rank pieces (we also refer
the reader to the remarks in \cite{AgaNW2011}).

However, in Theorem $1$ in the discussion by Ren and Zhou, a quite
stringent assumption requires that in some scaling regimes the true
low-rank component $L^\ast$ must vanish, that is,
$\|L^\ast\|_{\ell_\infty} \lesssim\sqrt{\frac{\log p}{n}}
\rightarrow
0$. The reason for this condition is effectively to ensure sign
consistency in recovering the sparse component. In a pure sparse model
selection problem (with no low-rank component in the population), the
deviation away from the sparse component is given only by noise due to
finite samples and this deviation is on the order of $\sqrt{\frac
{\log
p}{n}}$ in the Gaussian setting---consequently, sparse model selection
via $\ell_1$ norm thresholding is sign-consistent when the minimum
magnitude nonzero entry in the true model is larger than
$\sqrt{\frac{\log p}{n}}$. In contrast, if the true model consists of
both a sparse component and a low-rank component, the total deviation
away from the sparse component in the finite sample regime is given by
both sample noise as well as the low-rank component. This seems to be
the reason for the stringent assumption on the vanishing of the
low-rank component in Theorem $1$ of Ren and Zhou.

More broadly, one of the motivations of Ren and Zhou in proposing and
analyzing their estimator is that it may be possible to weaken the
assumptions on the minimum magnitude nonzero entry $\theta$ of the true
sparse component $S^\ast$ and the minimum nonzero singular value
$\sigma$ of the true low-rank component $L^\ast$---whether this is
possible under less stringent assumptions on $L^\ast$ is an interesting
question, and we comment on this point in Section \ref{sec:rates} in
the more general context of potentially improving the rates in our
paper.

\subsection{Other proposals}


Giraud and Tsybakov propose two alternative estimators for decomposing
a concentration matrix into sparse and low-rank components. While our
approach (\ref{eq:sdp}) builds on the graphical lasso, their proposed\vadjust{\goodbreak}
approaches build on the Dantzig selector of Cand\`es and Tao
\cite{CanT2007} and the neighborhood selection approach of Meinshausen
and B\"{u}hlmann \cite{MeiB2006}. Several comments are in order here.

First, we note that the extension of neighborhood selection proposed by
Giraud and Tsybakov to deal with the low-rank component begins by
reformulating the neighborhood selection procedure to obtain a
``global'' estimator that simultaneously estimates all the
neighborhoods. This reformulation touches upon a fundamental aspect of
latent variable modeling. In many applications marginalization over
the latent variables typically induces correlations between most pairs
of observed variables---consequently, local procedures that learn
model structure one node at a time are ill-suited for latent variable
modeling. Stated differently, requiring that a matrix be sparse with
few nonzeros per row or column (e.g., expressing preference for
a graphical model with bounded degree) can be done by imposing
column-wise constraints. On the other hand, the constraint that a
matrix be low-rank is really a global constraint expressed by requiring
all minors of a certain size to vanish. Thus, any estimator for latent
variable modeling (in the absence of additional conditions on the
latent structure) must necessarily be global in nature.

Second, we believe that the reformulation based on the Dantzig selector
perhaps ought to have an additional constraint. Recall that the
Dantzig selector \cite{CanT2007} constrains the $\ell_\infty$ norm (the
dual norm of the $\ell_1$ norm) of the correlated residuals rather than
the $\ell_2$ norm of the residuals as in the lasso. As the dual norm
of our combined $\ell_1$/trace norm regularizer involves both an
$\ell_\infty$ norm and a spectral norm, the following constraint set
may be more appropriate in the Dantzig selector based reformulation of
Giraud and Tsybakov:
\[
\mathcal{G} = \{(S,L)  \dvtx  \|\Sigma_O^n(S+L) - I\|_{\ell_\infty} \leq
\gamma\lambda_n,   \|\Sigma_O^n(S+L) - I\|_{2} \leq\lambda_n \}.
\]

Finally, we note that the Dantzig selector of \cite{CanT2007} has the
property that its constraint set contains the lasso solution (with the
same choice of regularization/relaxation parameters). In contrast,
this property is not shared in general by the Dantzig selector
reformulation of Giraud and Tsybakov in relation to our regularized
maximum-likelihood estimator (\ref{eq:sdp}). It is unclear how one
might achieve this property via suitable convex constraints in a
Dantzig selector type reformulation of our estimator.

In sum, both of these alternative estimators deserve further study.

\section{Comments on rates}\label{sec:rates}


Several of the commentaries (Wainwright, Giraud and Tsybakov, Ren and
Zhou and Cand\`es and Soltanolkotabi) bring up the possibility of
improving the rates given in our paper. At the outset we believe that
$n \gtrsim p$ samples is inherent to the latent variable modeling
problem if spectral norm consistency is desired in the low-rank
component. This is to be expected since the spectral norm of the
deviation of a sample covariance from the underlying population
covariance is on the order of $\sqrt{\frac{p}{n}}$. However, some more
subtle issues remain.

Giraud and Tsybakov point out that one may be concerned purely with
estimation of the sparse component, and that the low-rank component may
be a ``nuisance'' parameter. While this is not appropriate in every
application, in problem domains where the conditional graphical model
structure of the observed variables is the main quantity of interest
one can imagine quantifying deviations in the low-rank component via
``weaker'' norms than the spectral norm---this may lead to consistent
estimates for the sparse component with $n \ll p$ samples. The
analysis in our paper does not rule out this possibility, and a more
careful investigation is needed to establish such results.

Ren and Zhou suggest that while $n \gtrsim p$ may be required for
consistent estimation, one may be able to weaken the assumptions on
$\theta$ and $\sigma$ (the minimum magnitude nonzero entry of the
sparse component and the minimum nonzero singular value of the low-rank
component, respectively). From the literature on sparse model
selection, a natural lower bound on the minimum magnitude nonzero entry
for consistent model selection is typically given by the size of the
noise measured in the $\ell_\infty$ norm (the dual of the $\ell_1$
regularizer). Building on this intuition, a natural lower bound that
one can expect in our setting on $\theta$ is $\frac{1}{\gamma}
\|\Sigma_O^n - \Sigma\|_{\ell_\infty}$, while a natural bound on
$\sigma$ would be $\|\Sigma_O^n - \Sigma\|_2$. The reason for this
suggestion is that
$\max\{\frac{\|S\|_{\ell_\infty}}{\gamma},\|L\|_2\}$
is the
dual norm of the regularizer used in our paper. Therefore, it may be
possible to only require $\theta\sim
\frac{1}{\gamma}\sqrt{\frac{\log p}{n}}$ and $\sigma\sim
\sqrt{\frac{p}{n}}$. However, one issue here is that the
$\ell_\infty$ norm bound kicks in when $n \gtrsim\log p$ with
probability approaching one polynomially fast, while the spectral norm
bound only kicks in when $n \geq p$ but holds with probability
approaching one exponentially fast. Thus (as also noted by Giraud and
Tsybakov), it may be possible that $n \gtrsim p$ is required for
overall consistent estimation, but that the assumption on $\theta$
could be weakened by only requiring that the probability of consistent
estimation approach one polynomially fast.

Cand\`es and Soltanolkotabi comment that it would be of interest to
establish an ``adaptivity'' property whereby if no low-rank component
were present, the number of samples required for consistent estimation
would boil down to just the rate for sparse graphical model selection,
that is, $n \sim\log p$. While such a feature would clearly be
desirable to establish for our estimator, one potential roadblock may
be that our estimator (\ref{eq:sdp}) ``searches'' over a larger classes
of models than just those given by sparse graphical models;
consequently, rejecting the hypothesis that the observed variables are
affected by any latent variables may require that $n \gg\log p$. This
question deserves further investigation and, as suggested by Cand\`es
and Soltanolkotabi, recent results on adaptivity could inform a more
refined analysis of our estimator.

Finally, Wainwright suggests the intriguing possibility that faster
rates may be possible if the low-rank component has additional
structure. For example, there may exist a sparse factorization of the
low-rank component due to special structure between the latent and
observed variables. In such settings the trace norm regularizer
applied to the low-rank component is not necessarily the tightest
convex penalty. In recent joint work by the authors and Recht
\cite{ChaRPW2010}, a general framework for constructing convex penalty
functions based on some desired structure is presented. The trace norm
penalty for inducing low-rank structure is motivated from the viewpoint
that a low-rank matrix is the sum of a small number of rank-one
matrices and, therefore, the norm induced by the convex hull of the
rank-one matrices (suitably scaled) is a natural convex regularizer as
this convex hull (the trace norm ball) carries precisely the kind of
facial structure required for inducing low-rank structure in matrices.
In this spirit, one can imagine constructing convex penalty functions
by taking the convex hull of \emph{sparse} rank-one matrices. While
this convex hull is in general intractable to represent, relaxations of
this set that are tighter than the trace norm ball could provide faster
rates than can be obtained by using the trace norm.

\section{Weakening of irrepresentability conditions}



%

Wainwright asks a number of insightful questions regarding the
potential for weakening our Fisher information based conditions. Giraud
and Tsybakov also bring up connections between our conditions and
irrepresentability conditions in previous papers on sparse model
selection \cite{RavWRY2008,ZhaY2006}.

In order to better understand if the Fisher information based
conditions stated in our paper are necessary, Wainwright raises the
question of obtaining a converse result by comparing to an oracle
method that directly minimizes the rank and the cardinality of the
support of the components. A~difficulty with this approach is that we
don't have a good handle on the set of matrices that are expressible as
the sum of a sparse matrix and a low-rank matrix. The properties of
this set remain poorly understood, and developing a better picture has
been the focus of research efforts in algebraic geometry
\cite{DrtSS2007} and in complexity theory \cite{Val1977}. Nonetheless,
a comparison to oracle estimators that have side information about the
support of the sparse component and the row/column spaces of the
low-rank component (in effect, side information about the tangent
spaces at the two components) appears to be more tractable. This is
closer to the viewpoint we have taken in our paper in which we consider
the question of identifiability of the components given information
about the underlying tangent spaces. Essentially, our Fisher
information conditions state that these tangent spaces must be
sufficiently transverse with respect to certain natural norms and in a
space in which the Fisher information is the underlying inner-product.
More generally, as also pointed out by Giraud and Tsybakov, the
necessity of Fisher information based conditions is an open question
even in the sparse graphical model selection setting considered in
\cite{RavWRY2008}. The experimental studies in~\cite{RavWRY2008}
describing comparisons to neighborhood selection in some simple cases
provide a good starting point.

Wainwright raises the broader question of consistent model selection
when transversality of the underlying tangent spaces does not hold.
One approach~\cite{AgaNW2011} is to quantify the level of
identifiability based on a ``spikiness'' condition. A more geometric
viewpoint may be that only those pieces of the sparse and low-rank
components that do not lie in the intersection of their underlying
tangent spaces are fundamentally identifiable and, therefore,
consistency should be quantified with respect to these identifiable
pieces.

Giraud and Tsybakov ask about the interpretability of our conditions
$\xi(T)$ and $\mu(\Omega)$. These quantities are geometric in nature
and relate to the tangent space conditions for identifiability. In
particular, they are closely related to (and bounded by) the
incoherence of the row/column spaces of the low-rank component and the
maximum number of nonzeros per row/column \cite{ChaSPW2011}. These
latter quantities have appeared in many papers on sparse graphical
model selection (e.g., \cite{MeiB2006,RavWRY2008}) as well as on
low-rank matrix completion \cite{CanR2009}, and computing them is
straightforward. In our previous work on matrix decomposition
\cite{ChaSPW2011}, we note that these quantities are bounded for
natural random families of sparse and low-rank matrices based on
results in \cite{CanR2009}.

\section{Experimental issues and applications}\label{sec:exp}


%

Lauritzen and Meinshausen as well as Giraud and Tsybakov raise several
points about the choice of the regularization parameters. Choosing
these parameters in a data-driven manner (e.g., using the
methods described in \cite{MeiB2010}) is clearly desirable. We do wish
to emphasize that the sensitivities of the solution with respect to the
parameters $\lambda_n$ and $\gamma$ are qualitatively different. As
described in our main theorem and in our experimental section, the
solution of our estimator (\ref{eq:sdp}) is stable for a range of
values of $\gamma$ (see also \cite{ChaSPW2011})---this point is
observed by Yuan as well in his experiments. Further, the choice of
$\gamma$ ideally should not depend on $n$, while the choice of
$\lambda_n$ clearly should.

On a different point regarding experimental results, Giraud and
Tsybakov suggest at the end of their discussion that latent variable
models don't seem to provide significantly more expressive power than a
sparse graphical model. In contrast, Yuan's synthetic experiment seems
to provide compelling evidence that our approach (\ref{eq:sdp})
provides better performance relative to models learned by the graphical
lasso. The reason for these different observations may be tied to the
manner in which their synthetic models were generated. Specifically,
latent variable model selection using (\ref{eq:sdp}) is likely to be
most useful when the latent variables affect many observed variables
upon marginalization (e.g., latent variables are connected to many
observed variables), while the conditional graphical model among the
observed variables conditioned on the latent variables is sparse and
has bounded degree. This intuition is based on the theoretical
analysis in our paper and is also the setting considered in the
experiment in Yuan's discussion (as well as in the synthetic
experiments in our paper). On the other hand, the experimental setup
followed by Giraud and Tsybakov seems to generate a graphical model
with large maximum degree and low average degree, and randomly selects
a subset of the variables as latent variables. It is not clear if
these latent variables are the ones with large degree, which may
explain their remarks.

Finally, we note that sparse and low-rank matrix decomposition is
relevant in applications beyond the one described in our paper. As
observed by Lauritzen and Meinshausen, a natural matrix decomposition
problem involving \emph{covariance} matrices may arise if one considers
directed latent variable models in the spirit of factor analysis. In
such a context the covariance matrix may be expressed as the sum of a
low-rank matrix and a sparse (rather than just diagonal) matrix,
corresponding to the setting in which the distribution of the observed
variables conditioned on the latent variables is given by a sparse
covariance matrix. More broadly, similar matrix decomposition problems
arise in domains beyond statistical estimation such as optical system
decomposition, matrix rigidity and system identification in control
\cite{ChaSPW2011}, as well as others as noted by Cand\`es and
Soltanolkotabi.

\section{Future questions}\label{sec:future}

%
%
%
%

Our paper and the subsequent discussions raise a number of research and
computational challenges in latent variable modeling that we wish to
highlight briefly.

\subsection{Convex optimization in R} As mentioned by Lauritzen and
Meinshausen, R remains the software of choice for practitioners in
statistics. However, some of the recent advances in high-dimensional
statistical estimation have been driven by sophisticated convex
optimization based procedures that are typically prototyped using
packages such as SDPT3 \cite{TohTT} and others in Matlab and Python. It
would be of general interest to develop packages to invoke SDPT3
routines directly from R.


\subsection{Sparse/low-rank decomposition as infimal convolution} Given
a matrix $M \succ0$, consider the following function:
%
\begin{equation}\label{eq:decomp}
\|M\|_{S/L,\gamma} = \min_{S,L}   \gamma\|S\|_{\ell_1} + \tr(L),
\qquad\mbox{s.t. }   M = S-L,   L \succeq0.
\end{equation}
It is clear that $\|\cdot\|_{S/L,\gamma}$ is a norm, and it can be
viewed as the infimal convolution \cite{Roc1996} of the (scaled)
$\ell_1$ norm and the trace norm. In essence, it is a norm whose
minimization induces matrices expressible as the sum of sparse and
low-rank components (see also the atomic norm viewpoint of~\cite{ChaRPW2010}).\vadjust{\goodbreak} We could then effectively restate (\ref{eq:sdp})
as
\[
\hat{M}_n= \arg\min_{M \succ0}  -\lh(M; \Sigma^n_O)    + \lambda_n
\|M\|_{S/L,\gamma}
\]
and then decompose $\hat{M}_n$ by solving (\ref{eq:decomp}). This
two-step approach suggests the possibility of decoupling the
decomposition problem from the conditions fundamentally required for
consistency via regularized maximum-likelihood, as the latter only
ought to depend on the composite norm $\|\cdot\|_{S/L,\gamma}$. This
decoupling also highlights the different roles played by the parameters
$\lambda_n$ and $\gamma$ (as discussed in Section \ref{sec:exp}). More
broadly, such an approach may be useful as one analyzes general
regularizers, for example, convex penalties other than the trace norm
as described in Section \ref{sec:rates}.

\subsection{Non-Gaussian latent variable modeling}
As described in our paper and as raised by Wainwright, latent variable
modeling with non-Gaussian variables is of interest in many
applications. Both the computational and algebraic aspects present
major challenges in this setting. Specifically, the secant varieties
arising due to marginalization in non-Gaussian models (e.g., in models
with categorical variables) are poorly understood, and computing the
likelihood is also intractable. An approach based on matrix
decomposition as described in our paper may be appropriate, although
one would have to quantify the effects of the Gaussianity assumption.



\printaddresses

\end{document}